.
.
.
\font\sets=msbm10.
\font\stampatello=cmcsc10.
\font\symbols=msam10.

\def\1{{\bf 1}}
\def\sgn{{\rm sgn}}
\def\symsum{\sum_{|n-x|\le h}}

\def\starsum{\mathop{\enspace{\sum}^{\ast}}}
\def\dashsum{\mathop{\enspace{\sum}'}}

\def\square{\hbox{\vrule\vbox{\hrule\phantom{s}\hrule}\vrule}}
\def\defineq{\buildrel{def}\over{=}}
\def\defin{\buildrel{def}\over{\Longleftrightarrow}}

\def\multiplesum{\mathop{\sum \enspace \cdots \enspace \sum}}

\def\supporto{{\rm supp}\,}
\def\C{\hbox{\sets C}}
\def\N{\hbox{\sets N}}

\def\R{\hbox{\sets R}}
\def\Z{\hbox{\sets Z}}

\def\EssBdd{\hbox{\symbols n}\,}

\par
\centerline{\bf ON SOME LOWER BOUNDS}
\centerline{\bf OF SOME SYMMETRY INTEGRALS}
\bigskip
\centerline{by G.Coppola}
\bigskip
{
\font\eightrm=cmr8
\eightrm {
\par
{\bf Abstract.} We study the \lq \lq symmetry integral\rq \rq, \thinspace say $I_f$, of some arithmetic functions $f:\N \rightarrow \R$; we obtain from lower bounds of $I_f$ (for a large class of arithmetic functions $f$) lower bounds for the \lq \lq Selberg integral\rq \rq \thinspace of $f$, say $J_f$ (both these integrals give informations about $f$ in almost all the short intervals $[x-h,x+h]$, when $N\le x\le 2N$). In particular, when \thinspace $f=d_k$, the divisor function (having Dirichlet series \thinspace $\zeta^k$, with \thinspace $\zeta$ \thinspace the Riemann zeta function), where $k\ge 3$ is integer, we give lower bounds for the Selberg integrals, say \thinspace $J_k=J_{d_k}$, of the \thinspace $d_k$. We apply elementary methods (Cauchy inequality to get Large Sieve type bounds) in order to give $I_f$ lower bounds.
}
\footnote{}{\par \noindent {\it Mathematics Subject Classification} $(2010) : 11N37, 11N25.$}
}
\bigskip
\par
\noindent \centerline{\bf 1. Introduction and statement of the results.}
\smallskip
\par
\noindent
We give lower bounds of {\stampatello symmetry integrals} (here $\sgn(0)\defineq 0$, $r\neq 0$ $\Rightarrow$ $\sgn(r)\defineq {r\over {|r|}}$)
$$
I_f(N,h)\defineq \int_{N}^{2N} \Big| \symsum \sgn(n-x)f(n)\Big|^2 dx
$$
\par
\noindent
for a large class of arithmetic functions $f:\N \rightarrow \R$. For a motivation to study $I_f$, see esp. [C]. 
\par
A related integral is the, say, {\stampatello Selberg integral}, defined as 
$$
J_f(N,h)\defineq \int_{N}^{2N} \Big| \sum_{x<n\le x+h}f(n)-M_f(x,h)\Big|^2 dx, 
$$
\par
\noindent
where the {\stampatello mean-value} $M_f(x,h)$ depends \lq \lq weakly\rq \rq \thinspace on $x$ and is expected to depend linearly on $h$ (esp., it's of the kind $h$ times a polynomial in $\log x$, see the following). It is a kind of \lq \lq {\stampatello main term}\rq \rq \thinspace of the sum {\stampatello in the \lq \lq short interval\rq \rq} \thinspace $[x,x+h]$ (i.e., $h=o(x)$); so, we may expect it to approximate ($x\in \N$, $x\to \infty$ here), when $f=g\ast \1$, (compare [C1]) 
$$
h\Big({1\over x}\sum_{n\le x}f(n)\Big)={h\over x}\sum_d g(d)\left[ {x\over d}\right]\approx h\sum_{d\le x} {{g(d)}\over d}. 
$$
\par
While the former integral measures the {\stampatello almost-all} (i.e., for all $N\le x\le 2N$, except $o(N)$ of them) symmetry (around $x$) of $f$ in the {\stampatello short} (since $h=o(x)$) {\stampatello interval} $[x-h,x+h]$, the Selberg integral gives an \lq \lq {\stampatello average value}\rq \rq \thinspace to $f$ in $[x,x+h]$, for {\stampatello a.a.} (abbrev. almost all, s.i. shortens short intervals) these {\stampatello s.i.}
\par
Actually, it is a matter of evidence that knowing the (average) values of $f$ into {\stampatello a.a.s.i.} gives immediate information about the relative symmetry of $f$; however, let's go into more precise details and let's give an {\stampatello explicit connection} between these two integrals: 
$$
I_f(N,h)=\int_{N}^{2N}\Big| \sum_{x<n\le x+h}f(n)-\sum_{x-h\le n<x}f(n)\Big|^2 dx 
\ll \int_{N}^{2N}\Big| M_f(x,h)-M_f(x-h,h)\Big|^2 dx + \int_{N}^{2N}|f(x)|^2 dx + 
$$
$$
+ \int_{N-h}^{2N-h}|f(x)|^2 dx + \int_{N}^{2N}\Big| \sum_{x<n\le x+h}f(n)-M_f(x,h)\Big|^2 dx + \int_{N}^{2N}\Big| \sum_{x-h<n\le x}f(n)-M_f(x-h,h)\Big|^2 dx; 
$$
\par
\noindent
and, using the \lq \lq {\stampatello modified Vinogradov notation}\rq \rq, i.e. (in general, $F:\N \rightarrow \C$, here)
$$
F(N,h)\EssBdd G(N,h) \enspace \defin \enspace \forall \varepsilon>0 \enspace \enspace |F(N,h)|\ll_{\varepsilon} N^{\varepsilon}G(N,h), 
$$
\par
\noindent
so to leave (arbitrarily) small powers, {\stampatello assuming now on $f$ essentially bounded}, i.e. $f(n)\EssBdd 1$, we have 
$$
I_f(N,h)\ll_{\varepsilon} J_f(N,h) + \int_{N-h}^{2N-h}\Big|\sum_{x<n\le x+h}f(n)-M_f(x,h)\Big|^2 dx 
 + \int_{N}^{2N}\Big|M_f(x,h)-M_f(x-h,h)\Big|^2 dx + N^{1+\varepsilon} 
$$
$$
\ll_{\varepsilon} J_f(N,h) + \int_{N}^{2N}\Big| M_f(x,h)-M_f(x-h,h)\Big|^2 dx + N^{\varepsilon}(N+h^3), 
$$
\par
\noindent
where the last remainder comes from \lq \lq {\stampatello tails}\rq \rq, i.e. terms $\EssBdd h^3$ (see that $M_f(x,h)\EssBdd h$ is a consequence of the previous remarks on $M_f$ and $f\EssBdd 1$). We may assume, of course, that the difference $M_f(x,h)-M_f(x-h,h)$ is {\stampatello a.a.} small (i.e., its mean-square is \lq \lq small\rq \rq), due to the fact (compare the above, about $M_f$ choice) that $M_f$ is \lq \lq weakly\rq \rq \thinspace dependent on $x$ (like the case $M_k$ following, for $f=d_k$ with generating Dirichlet series $\zeta^k$). 
\medskip
\par
Then, ignoring these contributes together with the negligible $\EssBdd N+h^3$, we derive a lower bound of $J_f$, starting from a lower bound of $I_f$, here.
\medskip
\par
(We'll give a more precise calculation, following, for the more interesting cases $f=d_k$, see the above.)
\medskip
\par
We start simply remarking that the definition of \lq \lq {\stampatello mixed symmetry integrals}\rq \rq \thinspace (compare [C5]): 
$$
I_{f,f_1}(N,h)\defineq \int_{N}^{2N}\sum_{|n-x|\le h}\sgn(n-x)f(n)\sum_{|m-x|\le h}\sgn(m-x)f_1(m)dx 
$$
\par
\noindent
allows us to give a lower bound to \thinspace $I_f$, applying (expand the inner square), abbrev. $I_f$ for $I_f(N,h)$ \& similia, 
$$
0\le I_{f-f_1}=I_f-2I_{f,f_1}+I_{f_1}
$$
\par
\noindent
to get ($\forall N,h$ which are feasible)
$$
I_f(N,h)\ge 2I_{f,f_1}(N,h)-I_{f_1}(N,h).
\leqno{(1)}
$$
\par
\noindent
Here $(1)$ is true $\forall f,f_1:\N \rightarrow \R$ ({\stampatello any couple of real arithmetic functions}). 
\medskip
\par
However, in order to give a non-trivial lower bound to $I_f$, we need $I_{f_1}$ to be \lq \lq {\it smaller}\rq \rq \thinspace than $2I_{f,f_1}$. (That's the reason why we will give our general lower bound for \lq \lq mixed\rq \rq \thinspace integrals, but not for \lq \lq pure\rq \rq \thinspace ones.)
\medskip
\par
It will turn out, from our general result (next Theorem, compare the Lemma at next section), that the choice $f=d_k$ (general $k-$divisor function) and $f_1=d$ (i.e., $k=2$, divisor function) gives non-trivial lower bounds for $d_k$ symmetry integral; then, previous connection implies lower bounds for its Selberg integral, $J_k:=J_{d_k}$. This lower bound for $J_k$ is (ignoring logarithms) of the same order of magnitude of the diagonal (compare [C4], where this order of magnitude is required as an upper bound, to treat $2k-$th moments of $\zeta$).

\bigskip

\par
In order to simplify the exposition, we need to compare our variables to our {\stampatello main variable}, i.e. $N\to \infty$, from the point of view of exponents, using, say, $L:=\log N$ (\lq \lq logarithmic scale\rq \rq) :
\medskip
\item{a)} $\theta:={{\log h}\over L}$ is the {\stampatello width} (not the length, that's $h$) of the short interval $[x,x+h]$ (say, also of $[x-h,x+h]$);
\smallskip
\item{b)} $\lambda:={{\log Q}\over L}$ is the {\stampatello level} (see $\S1$ in [C3]) of \thinspace \thinspace $f:\N \rightarrow \R$, $f=g\ast \1$, $g(q)=0$, $\forall q>Q$;
\smallskip
\item{c)} $\delta:={{\log D}\over L}$ is the \lq \lq auxiliary level\rq \rq \thinspace of our mixed integral (in the following Theorem).
\medskip

\par
We explicitly remark that {\stampatello any inequality} involving these quantities {\stampatello will be} implicitly assumed to be {\stampatello sharp} (compare $\S1$ of [C3]): esp., our width will always be positive (i.e., $\exists \varepsilon_0>0$, absolute, with $\theta>\varepsilon_0>0$.). 
\medskip
\par
Our methods are elementary, as we apply a kind of {\stampatello Large Sieve Inequality}, using the {\stampatello spacing property of Farey fractions} (see the Lemma at next section).
\medskip
\par
We indicate, as usual, the {\stampatello distance to integers} (of any $\alpha \in \R$) as \thinspace $\Vert \alpha \Vert:=\min_{n\in \Z}|\alpha-n|$. 
\medskip
\par
Our results are the following. 
\smallskip

\par
\noindent {\stampatello Theorem.} {\it Fix } {\stampatello width} $0<\theta<1/2$, {\stampatello level} $0<\lambda<1$ {\it and } \lq \lq {\stampatello auxiliary level}\rq \rq \enspace $\delta$, {\it with } $\theta<\delta<\lambda$ {\it and } $\delta+\lambda<1$. {\it Let } $N,h,D,Q\in \N$, {\it with } $h=[N^{\theta}]$, $D=[N^{\delta}]$, $Q=[N^{\lambda}]$. {\stampatello Assume } $g_1:\N \rightarrow \R$, $g:\N \rightarrow \R$ {\it supported (resp.) in } $[1,D]$, $[1,Q]$, {\it with both } $1\le g_1\EssBdd 1$ {\it and } $1\le g\EssBdd 1$; {\it set } $f_1:=g_1\ast \1$, $f:=g\ast \1$. {\it Then, defining the} {\stampatello Ramanujan coefficients} {\it of an} {\stampatello essentially bounded} {\it arithmetic function} $F:§\N \rightarrow \C$ {\it as}
$$
R_{\ell}(F)\defineq \sum_{{m=1}\atop {m\equiv 0(\ell)}}^{\infty}{{(F\ast \mu)(m)}\over m}={1\over {\ell}}\sum_{n=1}^{\infty}{{G(\ell n)}\over n}
\EssBdd {1\over {\ell}},
$$
\par
\noindent
{\it where, say,} \thinspace $G\defineq F\ast \mu \EssBdd 1$ \enspace {\stampatello has finite support} ({\stampatello so} {\it to ensure} {\stampatello absolute convergence}), {\it we have} 
$$
I_{f,f_1}(N,h)=2N\sum_{1<\ell\le D}\ell^2 R_{\ell}(f)R_{\ell}(f_1)\sum_{t|\ell}{{\mu(t)}\over {t^2}}\Big\Vert {h\over {\ell/t}}\Big\Vert + o(Nh); 
$$
\par
\noindent
{\it whence, with an absolute constant, }
$$
I_{f,f_1}(N,h)\gg N\sum_{1<\ell\le {D\over 2}}\sum_{d\le {D\over {\ell}}}{{g_1(\ell d)}\over d}\sum_{q\le {Q\over {\ell}}}{{g(\ell q)}\over q}\sum_{t|\ell}{{\mu(t)}\over {t^2}}\Big\Vert {h\over {\ell/t}}\Big\Vert. 
$$
\par
\noindent
{\it Furthermore, assuming also that } $g(\ell q)\ge g(q)$ \thinspace $\forall \ell \le Q$, $\forall q\le {Q\over {\ell}}$, \thinspace \thinspace {\it we get the} ({\it absolute}) {\it lower bound} 
$$
I_{f,f_1}(N,h)\gg N\Big(\sum_{q\le {Q\over D}}{{g(q)}\over q}\Big)\sum_{1<\ell\le {D\over 2}}\sum_{d\le {D\over {\ell}}}{{g_1(\ell d)}\over d}\sum_{t|\ell}{{\mu(t)}\over {t^2}}\Big\Vert {h\over {\ell/t}}\Big\Vert. 
$$

\medskip
\par
Our \lq \lq main\rq \rq \thinspace consequence is for the symmetry integral of $d_k$ and for its Selberg integral, $J_k$, in the:
\smallskip

\par
\noindent {\stampatello Corollary.} {\it Fix } $k\ge 3$ \thinspace {\stampatello integer}. {\it Let } $N,h\in \N$ \thinspace {\it give, say, } {\stampatello width} \enspace $\theta=\theta_k$, $0<\theta_k<1/k$. {\it Then} 
$$
I_{d_k}(N,h)\gg_k NhL^{k+1},
\quad
J_k(N,h)\gg_k NhL^{k+1}.
$$

\medskip

\par
We explicitly remark the coincidence that the width $<{1\over k}$ is the range of $h$ for which the $J_k$ upper bound of the kind above (but it's a lower one !) is required, in order to get the (highly!) non-trivial bound (in [C4]) of $\zeta^{2k}$ integral-mean.

\bigskip

\par
\noindent
The paper is organized as follows:
\medskip
\item{$\diamond$} in section $2$ we state and prove our Lemma (on a \lq \lq discrete mixed integral\rq \rq);
\smallskip
\item{$\diamond$} in section $3$ we apply the Lemma (and an asymptotic formula) to prove our Theorem;
\smallskip
\item{$\diamond$} last section contains the proof of the Corollary, with some comments and remarks.

\bigskip
\bigskip
\bigskip

\par
\noindent \centerline{\bf 2. Statement and Proof of the Lemma.}
\smallskip
\par
\noindent
Our Lemma, following, deals with \lq \lq {\stampatello mixed symmetry integrals}\rq \rq, defined above as: 
$$
I_{f,f_1}(N,h)=\int_{N}^{2N}\sum_{|n-x|\le h}\sgn(n-x)f(n)\sum_{|m-x|\le h}\sgn(m-x)f_1(m)dx, 
$$
\par
\noindent
where $f=g\ast \1$, $f_1=g_1\ast \1$; actually, the mean-square in the Lemma is a discrete one (a sum !), not an integral (like in the previous version of this paper). This is done in order to apply Lemma 2 in [C-S] (a kind of Large Sieve inequality, see its proof), dealing with {\stampatello Farey fractions} (i.e., $j/\ell$, $r/t$, with $(j,\ell)=1=(r,t)$, see the proof) and exploiting their \lq \lq {\stampatello well-spaced}\rq \rq \thinspace  property (compare $(\ast)$ in the proof of the Lemma).
\par
By the way, the first appearance of these (\lq \lq mixed\rq \rq) integrals is in [C5], where (from the Cauchy-Schwarz inequality) they have non-trivial bounds, whenever one of the two \lq \lq pure\rq \rq \thinspace (symmetry) integrals has one: 
$$
|I_{f,f_1}(N,h)|\le \int_{N}^{2N}\Big|\sum_{|n-x|\le h}\sgn(n-x)f(n)\Big| \thinspace \Big|\sum_{|m-x|\le h}\sgn(m-x)f_1(m)\Big|dx 
\le \sqrt{I_f(N,h)} \sqrt{I_{f_1}(N,h)} .
$$
\par
Furthermore, we recall that the proof of the Lemma we use from our Acta Arithmetica paper relies solely on the Cauchy inequality. Hence, the present  Lemma inherits the elementary character from that one. 
\par
In fact, it comes from the properties (see [C-S]) of the function
$$
\chi_q(x)\defineq \sum_{{|n-x|\le h}\atop {n\equiv 0(\!\!\bmod q)}}\sgn(n-x), 
$$
\par
\noindent
entering the game, since (when $g,g_1$ have supports $\supporto(g_1)\subset [1,D]$, $\supporto(g)\subset [1,Q]$, here) 
$$
\sum_{x\sim N}\Big(\sum_{|n-x|\le h}\sgn(n-x)\sum_{{q\le Q}\atop {q|n}}g(q)\Big)\Big(\sum_{|m-x|\le h}\sgn(m-x)\sum_{{d\le D}\atop {d|m}}g_1(d)\Big) 
=\sum_{x\sim N}\sum_{q\le Q}g(q)\chi_q(x)\sum_{d\le D}g_1(d)\chi_d(x). 
$$
\par
\noindent
This {\stampatello discrete mixed integral} is linked to $I_{f,f_1}(N,h)$, see Thm. proof ($\S3$). We treat the sum (and not the integral, as mistaken in v1, previous version !) of this double sum over these \lq \lq character-like\rq \rq \thinspace functions. However, the Lemma still holds for $I_{f,f_1}(N,h)$ (as stated in v1, but will be proved within Thm. proof in $\S3$). 
\medskip
\par
We can (with Ramanujan coefficients $R_{\ell}(f)$ defined in the Thm. above) state and show our
\smallskip
\par
\noindent {\stampatello Lemma.} {\it Let } $N,h\in \N$ {\it with } $h\to \infty$ {\it and } $h=o(N)$ {\it when } $N\to \infty$. {\stampatello Assume } $g,g_1:\N \rightarrow \C$ {\it with } $g_1(d)=0$ $\forall d>D$ {\it and } $g(q)=0$ $\forall q>Q$, {\it where } $1<D\le Q\ll N$. {\stampatello Then} 
$$
\sum_{x\sim N}\Big(\sum_{q\le Q}g(q)\chi_q(x)\sum_{d\le D}g_1(d)\chi_d(x)\Big)=2N\sum_{1<\ell\le D}\ell^2 R_{\ell}(f)R_{\ell}(f_1)\sum_{t|\ell}{{\mu(t)}\over {t^2}}\left\Vert {h\over {\ell/t}}\right\Vert + 
$$
$$
 + {\cal O}\left( DQL\sqrt{\sum_{1<t\le 2h}t^2 |R_t(f_1)|^2+h\sum_{2h<t\le D}t |R_t(f_1)|^2}
                      \sqrt{\sum_{1<\ell \le 2h}\ell^2 |R_{\ell}(f)|^2+h\sum_{2h<\ell \le Q}\ell |R_{\ell}(f)|^2}\right).
$$

\smallskip

\par
\noindent {\stampatello proof.}$\!$ Abbreviate $n\equiv a(q)$ for $n\equiv a(\bmod \, q)$ and, from {\stampatello additive characters orthogonality} [V], 
$$
\chi_q(x)=\sum_{{|r|\le h}\atop {r\equiv -x(q)}}\sgn(r)=\sum_{j<q}c_{j,q}e_q(jx)
$$
\par
get the {\stampatello Fourier coefficients} (of previous finite Fourier expansion), see [C-S], 
$$
c_{j,q}:={1\over q}\sum_{|r|\le h}\sgn(r)e_q(rj) 
\enspace \thinspace \hbox{\stampatello satisfying} \thinspace \enspace 
c_{dj',dq'}={1\over d}c_{j',q'}, \forall d,j',q'\in \N, 
\enspace \hbox{\stampatello whence}
$$
$$
\chi_q(x)=\sum_{{\ell|q}\atop {\ell>1}}{{\ell}\over q}\starsum_{j\le \ell}c_{j,\ell}e_{\ell}(jx), 
\enspace \hbox{\stampatello with} \enspace 
\sum_{j<q}\left|c_{j,q}\right|^2 = 2\left\Vert {h\over q}\right\Vert, 
\quad
\starsum_{j<\ell}\left|c_{j,\ell}\right|^2 = 2\sum_{t|\ell}{{\mu(t)}\over {t^2}}\left\Vert {{ht}\over {\ell}}\right\Vert.
$$
\medskip
\par
By the way, this last relation highlights: {\stampatello the sum above}, performed over $t|\ell$, {\stampatello is non-negative}. Then 
$$
\sum_{x\sim N}\sum_{d\le D}\sum_{q\le Q}g_1(d)g(q)\chi_d(x)\chi_q(x) = 
$$
$$
=\sum_{1<t\le D}\left(\sum_{d'\le {D\over t}}{{g_1(td')}\over {d'}}\right)
  \sum_{1<\ell \le Q}\left(\sum_{q'\le {Q\over {\ell}}}{{g(\ell q')}\over {q'}}\right)
    \starsum_{r(t)}\overline{c_{r,t}}
     \starsum_{j(\ell)}c_{j,\ell}\sum_{x\sim N}e(\alpha x) 
$$
\par
({\stampatello apply previous properties} of $\chi_q$ {\stampatello expansion}), with, say, \enspace $\alpha:={j\over {\ell}}-{r\over t}$; \thinspace  apply Lemma 2 [C-S], {\stampatello since}\enspace 
$$
\Vert \alpha \Vert \neq 0 \enspace \Rightarrow \enspace \sum_{x\sim N}e(\alpha x) = e(\alpha/2){{e(2N\alpha)-e(N\alpha)}\over {2i\sin \pi \alpha}}\ll {1\over {\Vert \alpha \Vert}} 
\leqno{(*)}
$$
\par
together with\enspace ${\displaystyle {j\over {\ell}}\neq {r\over t} }$ \thinspace $\Rightarrow $ \thinspace ${\displaystyle \left\Vert {j\over {\ell}}-{r\over t}\right \Vert \ge {1\over {\ell t}}\gg {1\over {DQ}} }$ \enspace $\forall t\le D$ $\forall \ell \le Q$ ({\stampatello recall they're Farey fractions}) 
\par
{\stampatello give ${1\over {DQ}}$ well-spaced (Farey) fractions and (isolating ${j\over {\ell}}={r\over t}$ $\Rightarrow $ $\ell=t$, i.e. the \lq \lq diagonal\rq \rq}) 
$$
\sum_{x\sim N}\Big(\sum_{q\le Q}g(q)\chi_q(x)\sum_{d\le D}g_1(d)\chi_d(x)\Big) 
=\sum_{1<\ell \le D}\left(\sum_{d\le {D\over {\ell}}}{{g_1(\ell d)}\over d}\right)
                                   \left(\sum_{q\le {Q\over {\ell}}}{{g(\ell q)}\over q}\right)
                                    \left(\starsum_{j<\ell}|c_{j,\ell}|^2\right)N + 
$$
$$
 + {\cal O}\left( DQL\sqrt{\sum_{1<t\le D}\left| \sum_{d\le {D\over t}}{{g_1(td)}\over d}\right|^2 \starsum_{r<t}|c_{r,t}|^2}
                      \sqrt{\sum_{1<\ell \le Q}\left| \sum_{q\le {Q\over {\ell}}}{{g(\ell q)}\over q}\right|^2 \starsum_{j<\ell}|c_{j,\ell}|^2}\right).
$$
\par
From the above property of $\chi_q$ we {\stampatello may} use : \enspace $\sum_{j<\ell}^{*}|c_{j,\ell}|^2 \ll \sum_{j<\ell}|c_{j,\ell}|^2 \ll \min(1,{h\over {\ell}})$\enspace to get
$$
\sum_{1<t\le D}\left| \sum_{d\le {D\over t}}{{g_1(td)}\over d}\right|^2 \starsum_{r<t}|c_{r,t}|^2 
\ll \sum_{1<t\le 2h}\left| \sum_{d\le {D\over t}}{{g_1(td)}\over d}\right|^2 
     +h\sum_{2h<t\le D}{1\over t}\left| \sum_{d\le {D\over t}}{{g_1(td)}\over d}\right|^2 
$$
\par
and 
$$
\sum_{1<\ell \le Q}\left| \sum_{q\le {Q\over {\ell}}}{{g(\ell q)}\over q}\right|^2 \starsum_{j<\ell}|c_{j,\ell}|^2 
\ll \sum_{1<\ell \le 2h}\left| \sum_{q\le {Q\over {\ell}}}{{g(\ell q)}\over q}\right|^2 
     +h\sum_{2h<\ell \le Q}{1\over {\ell}}\left| \sum_{q\le {Q\over {\ell}}}{{g(\ell q)}\over q}\right|^2, 
$$
\par
{\stampatello whence} the (remainders, i.e. the) {\stampatello off-diagonal terms are} 
$$
\ll DQL 
     \sqrt{\sum_{1<t\le 2h}\left| \sum_{d\le {D\over t}}{{g_1(td)}\over d}\right|^2 
            +h\sum_{2h<t\le D}{1\over t}\left| \sum_{d\le {D\over t}}{{g_1(td)}\over d}\right|^2} 
\times 
$$
$$
\times
      \sqrt{\sum_{1<\ell \le 2h}\left| \sum_{q\le {Q\over {\ell}}}{{g(\ell q)}\over q}\right|^2 
             +h\sum_{2h<\ell \le Q}{1\over {\ell}}\left| \sum_{q\le {Q\over {\ell}}}{{g(\ell q)}\over q}\right|^2}; 
$$
\par
and, using the definition of {\stampatello Ramanujan coefficients} (see Thm.), we get the desired estimate$.\enspace \square$

\bigskip
\bigskip
\bigskip

\par
\noindent \centerline{\bf 3. Proof of the Theorem.}
\smallskip
\par
\noindent {\stampatello proof}. First of all, we link $I_{f,f_1}(N,h)$ with the discrete mixed integral of the Lemma: define 
$$
S_f^{\pm}(x)\defineq \symsum f(n)\sgn(n-x), 
$$
\par
the {\stampatello symmetry sum} of the (real) arithmetic function $f$. Obviously, $f\EssBdd 1$ $\Rightarrow $ $S_f^{\pm}\EssBdd h$. Then 
$$
I_{f,f_1}(N,h)=\int_{N}^{2N}\left(S_f^{\pm}([x])-f([x])+f([x]-h)\right)\left(S_{f_1}^{\pm}([x])-f_1([x])+f_1([x]-h)\right)dx =
$$
$$
= \sum_{N\le x<2N}\left(S_f^{\pm}(x)-f(x)+f(x-h)\right)\left(S_{f_1}^{\pm}(x)-f_1(x)+f_1(x-h)\right) = 
$$
$$
= \sum_{x\sim N}\left(S_f^{\pm}(x)-f(x)+f(x-h)\right)\left(S_{f_1}^{\pm}(x)-f_1(x)+f_1(x-h)\right) 
 + {\cal O}_{\varepsilon}\left(N^{\varepsilon}h^2\right). 
$$
\par
Here the $x$ is intended both real (in $\int$) and natural (in $\sum$); but the integral doesn't see the $x\in \N$.
\par
Due to the hypothesis $\theta<1/2$ ($\Rightarrow$ $\theta<1$), this error term is $o(Nh)$. Now, this sum is 
$$
\sum_{x\sim N}\left(\sum_{q\le Q}g(q)\chi'_q(x)\sum_{d\le D}g_1(d)\chi'_d(x)\right), 
$$
\par
which is not treated in the Lemma, because here ($x\in \N$ and $\1_{\wp}=1$ if $\wp$ is true, $0$ otherwise) : 
$$
\chi'_q(x)\defineq \dashsum_{{|n-x|\le h}\atop {n\equiv 0(\!\!\bmod q)}}\sgn(n-x) = \sum_{{x<n\le x+h}\atop {n\equiv 0(\!\!\bmod q)}}1-\sum_{{x-h<n\le x}\atop {n\equiv 0(\!\!\bmod q)}}1 = \chi_q(x)-\1_{q|x}+\1_{q|x-h}, 
$$
\par
i.e. the dash means that $n=x$ is counted with \lq \lq $-$\rq \rq \thinspace sign and $n=x-h$ is not counted. If we consider 
$$
\sum_{q\le Q}g(q)\chi'_q(x)-\sum_{q\le Q}g(q)\chi_q(x)=\sum_{q|x-h,q\le Q}g(q)-\sum_{q|x,q\le Q}g(q)\EssBdd d(x-h)+d(x)\EssBdd 1, 
$$
\par
we have that the present mean-square and the one in the Lemma differ by $\EssBdd Nh$. This is not negligible. 
\par
However, the same proof of the Lemma, applied to 
$$
\chi'_q(x) \enspace \hbox{\stampatello instead \thinspace of} \enspace \chi_q(x), \enspace \hbox{\stampatello with } c'_{j,q} \enspace \hbox{\stampatello instead \thinspace of} \enspace c_{j,q}, 
$$
\par
i.e. giving again the (finite) Fourier expansion, but with, say, the {\stampatello Fourier coefficients} 
$$
c'_{j,q}:={1\over q}\dashsum_{|r|\le h}\sgn(r)e_q(rj) 
$$
\par
(the dash takes $r=0$ with \lq \lq $-$\rq \rq \thinspace and doesn't count $r=-h$), we may repeat Lemma proof verbatim to 
$$
\chi'_q(x)=\sum_{{\ell|q}\atop {\ell>1}}{{\ell}\over q}\starsum_{j\le \ell}c'_{j,\ell}e_{\ell}(jx), 
\enspace \hbox{\stampatello with} \enspace 
\sum_{j<q}\left|c'_{j,q}\right|^2 = 2\left\Vert {h\over q}\right\Vert, 
\quad
\starsum_{j<\ell}\left|c'_{j,\ell}\right|^2 = 2\sum_{t|\ell}{{\mu(t)}\over {t^2}}\left\Vert {{ht}\over {\ell}}\right\Vert, 
$$
\par
getting (see the above; by the way, this recovers the Lemma in version v1) 
$$
I_{f,f_1}(N,h)=2N\sum_{1<\ell\le D}\ell^2 R_{\ell}(f)R_{\ell}(f_1)\sum_{t|\ell}{{\mu(t)}\over {t^2}}\left\Vert {h\over {\ell/t}}\right\Vert + o(Nh) + 
$$
$$
+ {\cal O}\left( DQL\sqrt{\sum_{1<t\le 2h}t^2 |R_t(f_1)|^2+h\sum_{2h<t\le D}t |R_t(f_1)|^2}
                      \sqrt{\sum_{1<\ell \le 2h}\ell^2 |R_{\ell}(f)|^2+h\sum_{2h<\ell \le Q}\ell |R_{\ell}(f)|^2}\right).
$$
\par
This holds for $\theta<1/2$ (as we use it for the Thm.), but is true in the Lemma hypotheses, joining $\theta<1$. 
\par
An immediate application of this gives Thm. equation, using $R_{\ell}(f),R_{\ell}(f_1)\EssBdd {1\over {\ell}}$ above, since $\delta+\lambda<1$. 
\par
Then, due to : $g_1,g\ge 1$, 
$$
2N\sum_{1<\ell\le D}\ell^2 R_{\ell}(f)R_{\ell}(f_1)\sum_{t|\ell}{{\mu(t)}\over {t^2}}\Big\Vert {h\over {\ell/t}}\Big\Vert 
= 2N\sum_{1<\ell\le D}\sum_{d\le {D\over {\ell}}}{{g_1(\ell d)}\over d}\sum_{q\le {Q\over {\ell}}}{{g(\ell q)}\over q}\sum_{t|\ell}{{\mu(t)}\over {t^2}}\Big\Vert {h\over {\ell/t}}\Big\Vert \gg 
$$
$$
\gg N\left(\sum_{q\le {Q\over D}}{1\over q}\right)\sum_{1<\ell\le D}\sum_{t|\ell}{{\mu(t)}\over {t^2}}\Big\Vert {h\over {\ell/t}}\Big\Vert \gg N\left(\sum_{q\le {Q\over D}}{1\over q}\right)\sum_{t<{D\over {2h}}}{{\mu(t)}\over {t^2}}\left(\sum_{1<n\le 2h}\Big\Vert {h\over n}\Big\Vert + h\sum_{2h<n\le {D\over t}}{1\over n}\right) \gg 
$$
$$
\gg Nh\log{Q\over D}\sum_{t<{D\over {2h}}}{{\mu(t)}\over {t^2}}\sum_{2h<n\le {D\over t}}{1\over n}
\gg Nh\log{Q\over D}\sum_{2h<n\le D}{1\over n}\sum_{t\le {D\over n}}{{\mu(t)}\over {t^2}} 
\gg Nh\log{Q\over D}\log{D\over {2h}}\gg NhL^2,  
$$
\par
where we used $\log{Q\over D}\gg L$ (from \thinspace $\delta<\lambda$) and $\log {D\over {2h}}\gg L$ (from \thinspace $\delta>\theta$) in the well-known (see, esp., [T]): 
$$
\sum_{t\le T}{{\mu(t)}\over {t^2}} = {1\over {\zeta(2)}} + {\cal O}\left( {1\over T}\right).
$$
\par
Hence, the main term is $\gg Nh$, whence $o(Nh)$ can be neglected (with an absolute constant in the $\gg$): 
$$
2N\sum_{1<\ell\le D}\ell^2 R_{\ell}(f)R_{\ell}(f_1)\sum_{t|\ell}{{\mu(t)}\over {t^2}}\Big\Vert {h\over {\ell/t}}\Big\Vert 
\gg N\sum_{1<\ell\le {D\over 2}}\sum_{d\le {D\over {\ell}}}{{g_1(\ell d)}\over d}\sum_{q\le {Q\over {\ell}}}{{g(\ell q)}\over q}\sum_{t|\ell}{{\mu(t)}\over {t^2}}\Big\Vert {h\over {\ell/t}}\Big\Vert.
$$
\par
Recall, always in our calculations, that \thinspace $g_1,g\ge 1$ \thinspace and the sum over \thinspace $t|\ell$ \thinspace is \thinspace $\ge 0$, see Lemma proof.
\par
For the same reasons, the additional hypothesis on $g$ gives at once \enspace $I_{f,f_1}(N,h)$ \enspace lower bound: 
$$
\thinspace \enspace \thinspace 
N\sum_{1<\ell\le {D\over 2}}\sum_{d\le {D\over {\ell}}}{{g_1(\ell d)}\over d}\sum_{q\le {Q\over {\ell}}}{{g(\ell q)}\over q}\sum_{t|\ell}{{\mu(t)}\over {t^2}}\Big\Vert {h\over {\ell/t}}\Big\Vert \gg N\Big(\sum_{q\le {Q\over D}}{{g(q)}\over q}\Big)\sum_{1<\ell\le {D\over 2}}\sum_{d\le {D\over {\ell}}}{{g_1(\ell d)}\over d}\sum_{t|\ell}{{\mu(t)}\over {t^2}}\Big\Vert {h\over {\ell/t}}\Big\Vert.\enspace \square 
$$

\bigskip
\bigskip
\bigskip

\par
\noindent \centerline{\bf 4. Proof of the Corollary. Remarks and comments.}
\smallskip
\par
\noindent {\stampatello proof}. We recall the definition of {\stampatello symmetry sum} for $f$ (see Thm. proof)  
$$
S_f^{\pm}(x)=\symsum \sgn(n-x)f(n)
$$
\par
\noindent
and, in particular, for $f=d_k$ (the $k-$divisor function, generated by $\zeta^k$), we write 
$$
S_k^{\pm}(x)\defineq \symsum d_k(n)\sgn(n-x) 
= \multiplesum_{{d_1 \enspace \qquad \enspace d_k}\atop {|d_1\cdots d_k-x|\le h}}\sgn(d_1\cdots d_k-x)
$$
\par
\noindent
where, considering that (in our symmetry sum) \thinspace 
\medskip
\par
\noindent
here \thinspace $x\ge N$ $\Rightarrow $ $x-h\ge N-h$ $\Rightarrow $ {\stampatello at least one of $d_1, \ldots , d_k$ has to be} \thinspace $d_j\ge (N-h)^{1/k}$, do the following: 
$$
\hbox{\stampatello let's\thinspace call}\enspace \Sigma_0 \enspace \hbox{\stampatello the\thinspace part\thinspace of}\enspace S_k^{\pm} \enspace \hbox{\stampatello in\thinspace which}\enspace d_1\ge (N-h)^{1/k}; \enspace \hbox{\stampatello remains}\enspace S_k^{\pm}-\Sigma_0, \thinspace \hbox{\stampatello in\thinspace which}\enspace d_2\ge (N-h)^{1/k}, 
$$
$$
\hbox{\stampatello let's\thinspace call \thinspace it}\enspace \Sigma_1; \enspace \hbox{\stampatello remains}\enspace S_k^{\pm}-\Sigma_0-\Sigma_1, \enspace \hbox{\stampatello in\thinspace which}\enspace d_3\ge (N-h)^{1/k}, \thinspace \hbox{\stampatello say}\enspace \Sigma_2, \enspace \hbox{\stampatello and\thinspace so\thinspace on}.
$$
\medskip
\par
\noindent
Since in \thinspace $S_k^{\pm}$ \thinspace {\stampatello at least one} of \thinspace $d_1, \ldots , d_k$ \thinspace {\stampatello has to be} \thinspace $\ge (N-h)^{1/k}$, {\stampatello we get} 
$$
S_k^{\pm}(x)=\sum_{q\le {{x+h}\over {(N-h)^{1/k}}}}d_{k-1}(q)\sum_{{\left| m-{x\over q}\right|\le {h\over q}}\atop {m\ge (N-h)^{1/k}}}\sgn\left( m-{x\over q}\right) + \sum_{q\le {{x+h}\over {(N-h)^{1/k}}}}d_{k-1}^{(1)}(q)\sum_{{\left| m-{x\over q}\right|\le {h\over q}}\atop {m\ge (N-h)^{1/k}}}\sgn\left( m-{x\over q}\right) + 
$$
$$
+ \cdots + \sum_{q\le {{x+h}\over {(N-h)^{1/k}}}}d_{k-1}^{(k-1)}(q)\sum_{{\left| m-{x\over q}\right|\le {h\over q}}\atop {m\ge (N-h)^{1/k}}}\sgn\left( m-{x\over q}\right), 
$$
\par
\noindent
where \enspace $d_{k-1}^{(0)}(q):=d_{k-1}(q)$ \enspace {\stampatello has restrictions on} $0$ {\stampatello factors}, \enspace 
$$
d_{k-1}^{(1)}(q):=\multiplesum_{{d_1 \enspace \qquad \enspace d_{k-1}}\atop {{d_1 \cdots d_{k-1} = q}\atop {d_1<(N-h)^{1/k}}}}1
$$
\par
\noindent
{\stampatello has on} $1$ {\stampatello factor},\enspace and \enspace $\forall j\le k-1$, 
$$
d_{k-1}^{(j)}(q):=\multiplesum_{{d_1 \enspace \qquad \enspace d_{k-1}}\atop {{d_1 \cdots d_{k-1} = q}\atop {d_1,\ldots,d_j<(N-h)^{1/k}}}}1
$$
\par
\noindent
{\stampatello has} $j$ {\stampatello factors with restrictions (which are independent of $x$ !)}. 
\par
\noindent
Hence, {\stampatello calling}\enspace $g(q):=\sum_{j=0}^{k-1}d_{k-1}^{(j)}(q)$, ({\stampatello depends} on \thinspace $(N-h)^{1/k}$, {\stampatello too}), {\stampatello get} \enspace $1\le g(q)\le kd_{k-1}(q)\EssBdd_k 1$. 
\par
\noindent
We {\stampatello obtain} immediately that the symmetry sum \thinspace $S_k^{\pm}(x)$ \thinspace equals (we'll ignore the constants $k-$dependence) 
$$
\sum_{q\le {{x+h}\over {(N-h)^{1/k}}}}g(q)\sum_{{\left| m-{x\over q}\right|\le {h\over q}}\atop {m\ge (N-h)^{1/k}}}\sgn\left( m-{x\over q}\right) 
= \sum_{q\le {{x-h}\over {(N-h)^{1/k}}}}g(q)\chi_q(x) + {\cal O}_{\varepsilon}\left( \sum_{{{x-h}\over {(N-h)^{1/k}}}<q\le {{x+h}\over {(N-h)^{1/k}}}}\sum_{\left| m-{x\over q}\right|\le {h\over q}}x^{\varepsilon}\right).
$$
\par
\noindent
In these remainders, $q>{{x-h}\over {(N-h)^{1/k}}}\gg N^{1-1/k}$ (as \thinspace $N\le x\le 2N$ \thinspace in the integral) gives (from our hypotheses \thinspace $k>2$ \thinspace and \thinspace $\theta<1/k$) that \thinspace $h=o(q)$, whence the interval \thinspace $[{{x-h}\over q},{{x+h}\over q}]$ \thinspace contains (at most) one integer $m$ (the $m-$sum is \lq \lq {\stampatello sporadic}\rq \rq) and this, in turn, implies that the remainders are : 
$$
\EssBdd \sum_{{{x-h}\over {(N-h)^{1/k}}}<q\le {{x+h}\over {(N-h)^{1/k}}}}\sum_{\left| m-{x\over q}\right|\le {h\over q}}1
\EssBdd 1,
$$
\par
\noindent
since the $q-$sum, too, contains at most one integer ({\stampatello sporadicity} from: $\theta<1/k$).
\par
However, from Cauchy-Schwarz inequality, this contributes $\EssBdd Nh$, giving \lq \lq interference\rq \rq \thinspace with the lower bound (of the same order of magnitude, say {\stampatello diagonal-}like). We need a slight improvement on this bound for the remainder; this is done estimating its mean-square (recall, we're to find a lower bound for its $N\le x\le 2N$ integral!): bounding $S_k^{\pm}(x)\EssBdd h$ (trivially), the contribute in the integral due to these remainders becomes (apply the {\stampatello sporadicity argument} to the inner $q-$sum) 
$$
\EssBdd h\int_{N}^{2N}\sum_{{{x-h}\over {(N-h)^{1/k}}}<q\le {{x+h}\over {(N-h)^{1/k}}}}\sum_{\left| m-{x\over q}\right|\le {h\over q}}1\,dx
\EssBdd h\max_{{{N-h}\over {(N-h)^{1/k}}}<q\le {{2N+h}\over {(N-h)^{1/k}}}}\int_{N}^{2N}\sum_{\left| m-{x\over q}\right|\le {h\over q}}1\,dx 
$$
$$
\EssBdd h^2\max_{{{N-h}\over {(N-h)^{1/k}}}<q\le {{2N+h}\over {(N-h)^{1/k}}}}\sum_{{{N-h}\over q}\le m\le {{2N+h}\over q}}1 
\EssBdd Nh\Big({h\over {N^{1-1/k}}}\Big), 
$$
\par
\noindent
which is \thinspace $o(Nh)$, since (recall: \thinspace $k>2$) we have {\stampatello width} \thinspace $\theta<1/k<1-1/k$.
\par
Now on, we will ignore all of the \thinspace $o(Nh)$ \thinspace contributes to our integrals. 
\par
\noindent
Writing \thinspace \lq \lq $\sim$\rq \rq \thinspace to mean we're leaving (such) negligible remainders, we are left with 
$$
S_k^{\pm}(x)\sim \sum_{q\le {{x-h}\over {(N-h)^{1/k}}}}g(q)\chi_q(x) 
= \sum_{q\le Q}g(q)\chi_q(x) + \sum_{Q<q\le {{x-h}\over {(N-h)^{1/k}}}}g(q)\chi_q(x), 
$$
\par
\noindent
where we set \thinspace $Q:=(N-h)^{1-1/k}$; whence, we are enabled to say that \thinspace $\lambda:=1-1/k$ \thinspace is the {\stampatello level}. 
\par
In fact, the same arguments of our Lemma give the same estimates for non-diagonal terms in the case we have the further limitation \thinspace $q\le {{x-h}\over {(N-h)^{1/k}}}$, which depends on $x$, since $(\ast)$, in the Proof of the Lemma, holds whatever limitations hold on the summation interval; also, we get from the second sum a positive (better, non-negative) contribution, for our symmetry integral (say, \lq \lq on the diagonal\rq \rq). 
\medskip
\par
Hence, we will ignore the second sum (from a positivity argument, to be applied soon again). 
\medskip
\par
Finally, we may also ignore (see the remarks, following soon after) the parts inside $g$ having limitations on the factors (say, consider \thinspace 
$g(q):=d_{k-1}(q)$, here). In all, we are left, after applying the Theorem (with $\delta$ auxiliary level, $\theta_k<\delta<1/k$, $D:=[N^{\delta}]$ and $g_1=\1$; also, $d_{k-1}(\ell q)\ge d_{k-1}(q)$), to saying that $(1)$, together with the bound ([C2], compare [C-S]) \thinspace $I_{g_1\ast \1}(N,h)\ll NhL^3$ (thanks to exponent $3<k+1$, $\forall k>2$), gives 
$$
I_{d_k}(N,h)\gg_k N\left(\sum_{q\le {Q\over D}}{{d_{k-1}(q)}\over q}\right)\sum_{1<\ell \le {D\over 2}}\log{D\over {\ell}}\sum_{t|\ell}{{\mu(t)}\over {t^2}}\Big\Vert {h\over {\ell/t}}\Big\Vert \gg_k 
$$
$$
\gg_k N\left(\sum_{q\le {Q\over D}}{{d_{k-1}(q)}\over q}\right)\left(\sum_{1<n\le 2h}\Big\Vert {h\over n}\Big\Vert \sum_{t<{D\over {2h}}}{{\mu(t)}\over {t^2}}\log{D\over {nt}} + h\sum_{2h<n\le D}{1\over n}\sum_{t\le {D\over n}}{{\mu(t)}\over {t^2}}\log{D\over {nt}}\right) \gg_k 
$$
$$
\gg_k Nh\left(\sum_{q\le {Q\over D}}{{d_{k-1}(q)}\over q}\right)\sum_{2h<n\le D}{1\over n}\sum_{t\le {D\over n}}{{\mu(t)}\over {t^2}}\log{D\over {nt}} 
\gg_k NhL^2 \left(\sum_{q\le {Q\over D}}{{d_{k-1}(q)}\over q}\right), 
$$
\par
\noindent
as in Theorem proof ($0<\theta<\delta$), having used [T] 
$$
\sum_{2h<n\le D}{1\over n}\sum_{t\le {D\over n}}{{\mu(t)}\over {t^2}}\log{D\over {nt}} 
= \sum_{2h<n\le D}{1\over n}\log{D\over n}\left({1\over {\zeta(2)}}+{\cal O}\left({n\over D}\right)\right) - \sum_{2h<n\le D}{1\over n}\left(\sum_{t=1}^{\infty}{{\mu(t)\log t}\over {t^2}}+{\cal O}\left( {n\over D}L\right)\right), 
$$
\par
\noindent
together with partial summation [T] (compare [C-S] Corollary 1 calculations, p.199 on); hence 
$$
I_{d_k}(N,h)\gg_k NhL^{k+1},
$$
\par
\noindent
this last inequality coming from partial summation (see [D] or [T]) and (Lemma 1.1.2 of) Ch.1 of [L], as
$$
\sum_{n\le x}{{d_{k-1}(n)}\over n}=\multiplesum_{{n_1\enspace , \thinspace \ldots \thinspace ,\enspace n_{k-1}}\atop {n_1 \cdots n_{k-1} \le x}}{1\over {n_1 \cdots n_{k-1}}}\gg_k (\log x)^{k-1}.
$$
\par
\noindent
As regards the lower bound for Selberg integral $J_k$, we apply previous connection, with (in [C4] details) 
$$
M_k(x,h):=hP_{k-1}(\log x), 
$$
\par
\noindent
where $P_{k-1}(\log x)$ is a $k-1$ degree polynomial in $\log x$, whence 
$$
M_k(x,h)-M_k(x-h,h)=hM'_k(x-\alpha h,h)\EssBdd h^2/N, \qquad \forall x\in [N,2N]
$$
\par
\noindent
(from mean-value theorem, with $M'_k(x,h):={d\over {dx}}M_k(x,h)$, $0<\alpha<1$), whence (\lq \lq $\gg$\rq \rq \thinspace leaves \thinspace $o(Nh)$, here) 
$$
J_k(N,h)\gg I_{d_k}(N,h).\enspace \square
$$

\bigskip

\par
\noindent
We remark we \lq \lq wasted\rq \rq, in our (previous, version 1) lower bounds, \lq \lq {\it many}\rq \rq \thinspace terms in our previous analysis.
\medskip
\par
In fact, we felt that the limitation $k\ge 5$ was immaterial. 
\medskip
\par
Actually, the real improvement comes from the (previously) neglected terms of the Theorem, where the M\"obius function rendered more cumbersome our estimates (simplified by the hypothesis $g(\ell q)\ge g(q)$, here). 
\medskip
\par
Once again, we are postponing other eventual, further improvements to a future, forthcoming paper. 

\bigskip
\bigskip
\bigskip

\par
\noindent
\centerline{\bf References}
\smallskip
\item{\bf [C]} \thinspace Coppola, G.\thinspace - \thinspace {\sl On the symmetry of divisor sums functions in almost all short intervals} \thinspace - \thinspace Integers {\bf 4} (2004), A2, 9 pp. (electronic). $\underline{\tt MR\enspace 2005b\!:\!11153}$  
\smallskip
\item{\bf [C1]} \thinspace Coppola, G.\thinspace - \thinspace {\sl On the Correlations, Selberg integral and symmetry of sieve functions in short intervals} \thinspace - \thinspace http://arxiv.org/abs/0709.3648v3 (to appear on: Journal of Combinatorics and Number Theory)
\smallskip
\item{\bf [C2]} \thinspace Coppola, G.\thinspace - \thinspace {\sl On the Correlations, Selberg integral and symmetry of sieve functions in short intervals, II} \thinspace - \thinspace Int. J. Pure Appl. Math. {\bf 58.3}(2010), 281--298.
\smallskip
\item{\bf [C3]} \thinspace Coppola, G.\thinspace - \thinspace {\sl On the Correlations, Selberg integral and symmetry of sieve functions in short intervals, III} \thinspace - \thinspace http://arxiv.org/abs/1003.0302v1
\smallskip
\item{\bf [C4]} \thinspace Coppola, G.\thinspace - \thinspace {\sl On the Selberg integral of the $k-$divisor function and the $2k-$th moment of the Riemann zeta-function} \thinspace - \thinspace http://arxiv.org/abs/0907.5561v1 - to appear on Publ. Inst. Math., Nouv. Sér.
\smallskip
\item{\bf [C5]} \thinspace Coppola, G.\thinspace - \thinspace {\sl On the symmetry of arithmetical functions in almost all short intervals, V} \thinspace - \thinspace (electronic) http://arxiv.org/abs/0901.4738v2 
\smallskip
\item{\bf [C-S]} Coppola, G. and Salerno, S.\thinspace - \thinspace {\sl On the symmetry of the divisor function in almost all short intervals} \thinspace - \thinspace Acta Arith. {\bf 113} (2004), {\bf no.2}, 189--201. $\underline{\tt MR\enspace 2005a\!:\!11144}$
\smallskip
\item{\bf [D]} \thinspace Davenport, H.\thinspace - \thinspace {\sl Multiplicative Number Theory} \thinspace - \thinspace Third Edition, GTM 74, Springer, New York, 2000. $\underline{{\tt MR\enspace 2001f\!:\!11001}}$
\smallskip
\item{\bf [L]} \thinspace Linnik, Ju.V.\thinspace - \thinspace {\sl The Dispersion Method in Binary Additive Problems} \thinspace - \thinspace Translated by S. Schuur \thinspace - \thinspace American Mathematical Society, Providence, R.I. 1963. $\underline{\tt MR\enspace 29\# 5804}$
\smallskip
\item{\bf [T]} \thinspace Tenenbaum, G.\thinspace - \thinspace {\sl Introduction to Analytic and Probabilistic Number Theory} \thinspace - \thinspace Cambridge Studies in Advanced Mathematics, {\bf 46}, Cambridge University Press, 1995. $\underline{\tt MR\enspace 97e\!:\!11005b}$
\smallskip
\item{\bf [V]} \thinspace Vinogradov, I.M.\thinspace - \thinspace {\sl The Method of Trigonometrical Sums in the Theory of Numbers} - Interscience Publishers LTD, London, 1954. $\underline{{\tt MR \enspace 15,941b}}$

\medskip

\leftline{\tt Dr.Giovanni Coppola}
\leftline{\tt DIIMA - Universit\`a degli Studi di Salerno}
\leftline{\tt 84084 Fisciano (SA) - ITALY}
\leftline{\tt e-mail : gcoppola@diima.unisa.it}
\leftline{\tt e-page : www.giovannicoppola.name}

\bye